# Counting pairs of lattice paths by intersections

Ira Gessel[1], Wayne Goddard[2], Walter Shur, Herbert S. Wilf[3], and Lily Yen[4]

On an $r \times (n - r)$ lattice rectangle, we first consider walks that begin at the SW corner, proceed with unit steps in either of the directions E or N, and terminate at the NE corner of the rectangle. For each integer $k$ we ask for $N_k^{n,r}$, the number of *ordered* pairs of these walks that intersect in exactly $k$ points. The number of points in the intersection of two such walks is defined as the cardinality of the intersection of their two sets of vertices, excluding the initial and terminal vertices. The figure below shows a pair of such walks where $r = 9$, $n = 17$, and $k = 5$.

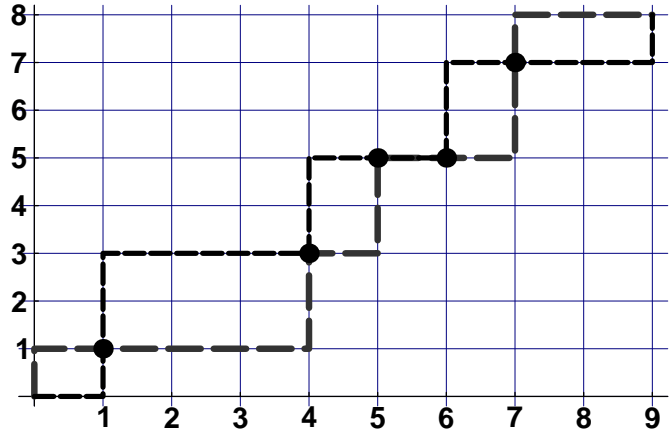

In sections 1 and 2 below, we will find two explicit formulas for the numbers $N_k^{n,r}$. Our formulas are

$$N_k^{n,r} = \frac{2(k+1)}{n-k-1} \sum_i \binom{k}{i} \binom{n-k+i-1}{r} \binom{n-i-1}{n-r}, \qquad (0 \le k \le n-2) \qquad (1a)$$

and

$$N_k^{n,r} = \frac{2(k+1)}{r} \sum_i (-1)^i \frac{\binom{k}{i}\binom{k-i}{i}\binom{n-i-2}{r-1}\binom{n-i-1}{r-i-1}}{\binom{n-i-2}{i}} \qquad (0 \le k \le n-2). \qquad (1b)$$

In section 2 we will prove the equality of (1a), (1b) in their common range of validity.

Next, we note that these formulas reveal the interesting fact that $N_1^{n,r} = 2N_0^{n,r}$, i.e., that *exactly twice as many pairs of walks have a single intersection as have no intersection*. Such a relationship clearly merits a bijective proof, and we will supply one in section 3 below.

---

[1] Supported by the National Science Foundation
[2] Supported by the Office of Naval Research
[3] Supported by the Office of Naval Research
[4] Supported by a Fellowship of the NSERC, Canada

In sections 4 and 5 we will discuss a number of related results. In the first of these sections we specify that the walkers end up at different points; in the second we ignore where the two walkers end up.

Finally, in section 6 we will discuss a variation in which we find the probability that two independent walkers on a given lattice rectangle do not meet, under a different hypothesis. In that situation the walkers start at the two points $(a, b+x+1)$ and $(a+x+1, b)$ in the first quadrant, and walk West or South at each step, except that when a walker reaches the $x$-axis (resp. the $y$-axis) then all future steps are constrained to be South (resp. West) until the origin is reached. We will find (Theorem 4 below) that if the probability $p(i,j)$ that a step from $(i,j)$ will go West depends only on $i+j$, *then the probability that the two walkers do not meet until they reach the origin is the same as the probability that a single (unconstrained) walker who starts at $(a, b+x+1)$ and takes $a+b+x$ steps, finishes at one of the points $(0,1), (-1,2), \ldots, (-x, 1+x)$.*

## 1. Derivation of the formula (1a)

We will use a generating function to approach the problem. If we sort out the pairs of walks that have $k$ intersections according to the $(m,q)$ of their last (i.e., most northeasterly) intersection, then we see the recurrence

$$N_k^{n,r} = \sum_{q,m} N_{k-1}^{m,q} N_0^{n-m, r-q}. \qquad (2)$$

Introduce the generating function $u_k(x,y) = \sum_{n,r} N_k^{n,r} x^n y^r$. Then (2) says simply that $u_k = u_{k-1} u_0$, for $k \geq 1$. Thus $u_k(x,y) = u_0(x,y)^{k+1}$ for $k = 0, 1, 2, \ldots$.

Now, $u_0$ is well known, having been calculated by Narayana [3]. (Indeed, the values $N_0^{n,r}/2$ are sometimes referred to as the Narayana numbers.) But it's more interesting to find $u_0$ by observing that the coefficient of $x^n y^r$ in $\sum_k u_k(x,y)$ is the total number of pairs of walks, since every pair has *some* number of intersections.

The number of all pairs of walks is $\binom{n}{r}^2$, so we have

$$1 + \sum_k u_k(x,y) = \sum_{k \geq 0} u_0(x,y)^k = \frac{1}{1 - u_0(x,y)}$$

$$= \sum_{n,r} \binom{n}{r}^2 x^n y^r$$

$$= \sum_{n \geq 0} x^n (y-1)^n P_n\left(\frac{y+1}{y-1}\right)$$

$$= \frac{1}{\sqrt{1 - 2x(y+1) + x^2(y-1)^2}},$$

where the $P_n$'s are the Legendre polynomials, and their classical generating function has been used. Thus
$$u_0(x,y) = 1 - \sqrt{1 - 2x(y+1) + x^2(y-1)^2}.$$

It follows that the number of pairs of walks that have exactly $k$ intersections is the coefficient of $x^n y^r$ in

$$u_k(x, y) = \left(1 - \sqrt{1 - 2x(y+1) + x^2(y-1)^2}\right)^{k+1}. \tag{3}$$

If we set $x = z$ and $y = y/z$, we can rephrase this conclusion as follows.

**Proposition.** *The number $N_k^{n,r}$ is the coefficient of $y^r z^{n-r}$ in $(y + z + 2f)^{k+1}$, where*

$$f = \frac{(1 - y - z) - \sqrt{(1 - y - z)^2 - 4yz}}{2}$$

*satisfies the equation $f = (y + f)(z + f)$.* ∎

To extract the coefficients of this generating function we use the Lagrange Inversion Formula for the solution of $f = xg(f)$ in the form (see, e.g., [5], eq. (5.1.2))

$$[x^n]\phi(f) = \frac{1}{n}[t^{n-1}]\phi'(t)g(t)^n. \tag{4}$$

First introduce an auxiliary variable $x$, and consider the equation $f = x(y + f)(z + f)$. Apply (4) using $\phi(t) = (y + z + 2t)^{k+1}$, and the formula (1a) follows.

## 2. Equality of (1a) and (1b)

This section is devoted to the proof of the equality of the two formulas (1a), (1b), in the range $0 \leq k \leq n - 2$. Note that we adopt the usual convention that any term with the factorial of a negative integer in its denominator is considered to be zero.

**Theorem 1.** *For $0 \leq k \leq n - 2$ we have*

$$\frac{k+1}{n-k-1} \sum_i \binom{k}{i} \binom{n-k+i-1}{r} \binom{n-i-1}{n-r}$$
$$= \frac{k+1}{r} \sum_j (-1)^j \frac{\binom{k}{j}\binom{k-j}{j}\binom{n-j-2}{r-1}\binom{n-j-1}{r-j-1}}{\binom{n-j-2}{j}}. \tag{5}$$

**Proof.** We obtain the two sides of (5) by evaluating in two ways the double sum

$$S = \sum_{j,l}(-1)^j \frac{(k+1)!\,(n-k-2)!\,(n-j-1)!}{r!\,(n-r)!\,j!\,l!\,(r-j-l-1)!\,(k-2j-l)!\,(n-k-r+j+l-1)!},$$

where the sum is over all nonnegative integers $j$ and $l$ with $j \leq n - 1$.

We shall need two forms of Vandermonde's theorem:

$$\sum_i \binom{a}{i}\binom{b}{m-i} = \binom{a+b}{m} \tag{6}$$

and
$$\sum_i (-1)^i \binom{a}{i}\binom{c-i}{m-i} = \binom{c-a}{m}. \tag{7}$$

First we sum on $l$. We have

$$S = \sum_{j=0}^{n-1}(-1)^j \frac{(k+1)!\,(n-j-1)!}{r!\,(n-r)!\,j!\,(k-2j)!} \sum_l \binom{k-2j}{l}\binom{n-k-2}{r-j-l-1}$$

$$= \sum_{j \le k/2}(-1)^j \frac{(k+1)!\,(n-j-1)!}{r!\,(n-r)!\,j!\,(k-2j)!}\binom{n-2j-2}{r-j-1} \quad \text{by (6),}$$

and this is easily seen to be equal to the right side of (5).

Next we set $l = i - j$ in $S$ and sum on $j$. This gives

$$S = \sum_{i,j}(-1)^j \frac{(k+1)!\,(n-k-2)!\,(n-j-1)!}{r!\,(n-r)!\,j!\,(i-j)!\,(r-i-1)!\,(k-i-j)!\,(n-k-r+i-1)!}$$

$$= \sum_i \frac{(k+1)!\,(n-k-2)!}{r!\,(n-r)!\,(r-i-1)!\,(n-k-r+i-1)!} \sum_{j=0}^{n-1}(-1)^j \frac{(n-j-1)!}{j!\,(i-j)!\,(k-i-j)!}$$

$$= \sum_i \frac{(k+1)!\,(n-k-2)!\,(n-i-1)!}{r!\,(n-r)!\,(r-i-1)!\,(n-k-r+i-1)!\,(k-i)!} \sum_{j=0}^{i}(-1)^j \binom{k-i}{j}\binom{n-j-1}{i-j}$$

$$= \sum_i \frac{(k+1)!\,(n-k-2)!\,(n-i-1)!}{r!\,(n-r)!\,(r-i-1)!\,(n-k-r+i-1)!\,(k-i)!}\binom{n-k+i-1}{i} \quad \text{by (7),}$$

and this is easily seen to be equal to the left side of (5).

It may be noted that the theorem can also be obtained from formula (1), p. 30 of [1] by taking the limit as $c \to \infty$ then setting $m = r - 1$, $a = -k$, $w = 1 - n$, and $b = n - k$.

Finally, we remark that Zeilberger's algorithm (see [6]) is, in principle, capable of verifying that the two sides of (5) satisfy the same recurrence relation. In fact, however, computers are not yet equal to the task of producing certificates for such two-variable summands as these in reasonable execution times.

## 3. A bijection

The formulas (1a,b) show that

$$N_0^{n,r} = \frac{2}{n-1}\binom{n-1}{r}\binom{n-1}{n-r},$$

and that $N_1^{n,r} = 2N_0^{n,r}$. Here we give a bijective proof of the latter assertion. For convenience, we define $s = n - r$ so that we deal with an $r \times s$ rectangle.

We will say that $(x, y)$ is an interior point of the rectangle if $0 < x < r$ and $0 < y < s$. The *distance at $x$* between two lattice paths $P_1$ and $P_2$ is
$$d_x(P_1, P_2) = \min\{|y_2 - y_1| : (x, y_1) \in P_1, (x, y_2) \in P_2\}.$$
We refer to eastbound and northbound edges of the paths as E edges and N edges, respectively. By the *shape* of a pair of lattice paths we mean the plane region enclosed between them.

**Definition of the map.** Now we define our map $\phi$ from pairs of paths with no intersections to pairs of pairs of paths with one intersection. Hence let $(P, Q)$ be a pair of paths that do not intersect. We can suppose that $P$ is north of $Q$. Suppose $d_x(P, Q) \geq 2$ for all $1 \leq x \leq r - 1$. Then $\phi$ maps $(P, Q)$ to $(P', Q)$ and $(P, Q')$, in which
  (i) $P'$ is obtained from $P$ as follows. First translate $P$ down by 1 unit, then delete its first N edge, and then concatenate an N edge to the last vertex of the new path. The pair $(P', Q)$ intersects at $(r, s-1)$ and only there.
  (ii) $Q'$ is obtained from $Q$ as follows. First translate $Q$ up by 1 unit, then delete its last N edge, then adjoin an N edge in-bound to its first vertex. The pair $(P, Q')$ intersects at $(0, 1)$, and only there.

Now suppose $P$ is north of $Q$ and $d_x(P, Q) = 1$ for some $x$, $1 \leq x \leq r - 1$. Let $x_0$ be the smallest such $x$, and let $y_0 = \max_y\{(x_0, y) \in Q\}$.

First suppose $(x_0, y_0) \neq (1, 0)$. Then lower by 1 unit the portion of $P$ from $(0, 0)$ to $(x_0, y_0 + 1)$, and move the first N edge of $P$ to join $(x_0, y_0)$ and $(x_0, y_0 + 1)$, to obtain the new path $P'$. The pair $(P', Q)$ intersects at $(x_0, y_0)$ and only there. To get another pair of paths that intersect at $(x_0, y_0)$, interchange the upper path with the lower path between $(x_0, y_0)$ and $(n, s)$, in $(P', Q)$.

Finally, suppose $(x_0, y_0) = (1, 0)$. Then we first produce $(P', Q)$ exactly as in the previous paragraph so that $(P', Q)$ intersects at $(1, 0)$ and only there. To produce a second pair that intersects at $(r - 1, s)$ and only there, the double-E edge from $(0, 0)$ to $(1, 0)$ in $(P', Q)$ is, in this case, moved to the northeast corner as another double-E edge. Then the resulting pair is translated 1 unit westward so the paths begin at $(0, 0)$, end at $(r, s)$, and intersect at $(r - 1, s)$ and only there.

**Invertibility of the map.** We partition the collection of all pairs of paths that intersect exactly once into groups of two as follows.
  (i) If $(P, Q)$ intersects at $(1, 0)$ then pair $(P, Q)$ with $(P', Q')$, where $(P', Q')$ intersects at $(r - 1, s)$ and the removal of the double-E edges from both pairs $(P, Q)$, $(P', Q')$ results in the same pair of nonintersecting paths on an $(r-1) \times s$ rectangle.
  (ii) If $(P, Q)$ intersects at $(0, 1)$ then pair $(P, Q)$ with $(P', Q')$, where $(P', Q')$ intersects at $(r, s - 1)$ and the removal of the double-N edges from both pairs $(P, Q)$, $(P', Q')$ results in the same pair of nonintersecting paths on an $r \times (s-1)$ rectangle.
  (iii) Suppose $(P, Q)$ intersects at an interior point. If $P$ is north of $Q$ from $(0, 0)$ to the intersection then pair $(P, Q)$ with $(P', Q')$ where the two paths have been interchanged from the intersection point to $(r, s)$, so $P'$ is always north of $Q'$.

Now we can define the inverse mapping $\psi = \phi^{-1}$. Given two pairs of paths $(P, Q)$, $(P', Q')$ from group (i) above. Then $\psi$ looks only at $(P, Q)$. It removes the intersection by lifting one of the double edges up by 1 unit in $(P, Q)$, and it moves the second edge of $P$, which is an E edge, to join $(0, 0)$ and $(0, 1)$.

For $(P, Q), (P', Q')$ from group (ii) above, $\psi$ looks at the pair $(P', Q')$ which intersects at $(r, s - 1)$. It lifts the upper path, except for its final N edge, by 1 unit, then moves that N edge to join $(0, 0)$ and $(0, 1)$.

Finally, for $(P, Q), (P', Q')$ from group (iii) above, $\psi$ looks at $(P, Q)$ such that $P$ is always north of $Q$, and suppose they meet at $(x_0, y_0)$. Then lift the portion of $P$ that precedes the intersection up by 1 unit, delete the N edge of $P$ from $(x_0, y_0)$ to $(x_0, y_0 + 1)$, and add an N edge to join $(0, 0)$ and $(0, 1)$. ∎

Narayana [4] showed that the value $N_0^{n,r}/2$ is equal to the number of plane trees with $n$ vertices and $r$ leaves. A proof of this using lattice paths is given in [2].

## 4. Termination at different endpoints

In this section we extend the formula for $N_k^{n,r}$ by considering pairs of walks where the two walkers start at the same point but end at different points. Say the walkers both start at $(0, 0)$, and the first walker terminates at $(r, n - r)$ and the second at $(s, n - s)$. Then, for $r < s$, let $M_{r,s}^{n,k}$ denote the number of (unordered) pairs of these walks that intersect in exactly $k$ points, not counting the starting point; for $r = s$, let $M_{r,r}^{n,k} = N_{k-1}^{n,r}$.

**Theorem 2.** *For $r \leq s$, the numbers $M_{r,s}^{n,k}$ are given by*

$$2\sum_t \sum_j (-1)^j \frac{(s - j - r + 1 + 2t)}{n - 1 - j - 2t} \binom{k}{2t + 1}\binom{k - 1 - 2t}{j}\binom{n - 1 - j - 2t}{s - j}\binom{n - 1 - j - 2t}{r - 1 - 2t}$$

$$+ \frac{s - r}{n - k} \sum_j \binom{k}{j}\binom{n - k}{r - j}\binom{n - k}{s - j}.$$

The second expression counts the number of pairs that intersect in each of the first $k$ steps; the other expression counts the remaining pairs. For $k = 0$ we recover the familiar:

$$M_{r,s}^{n,0} = \frac{s - r}{n} \binom{n}{r}\binom{n}{s}.$$

The proof is straightforward (at least conceptually) and we will omit the details. One first shows that the claimed formula is correct for $r = s$ by showing that the above expression is equivalent to the expression in (1b). Then for $r < s$, one easily obtains a recurrence relation for $M_{r,s}^{n,k}$ by considering the last and next-to-last steps of both walks. The remainder of the proof is to show that the above formula satisfies the recurrence and meets the proper boundary conditions.

## 5. Further remarks

In this section we discuss some variations of the original problem.

We start with the case where the walks are not restricted to a given lattice rectangle. More precisely, fix $n > 0$ and consider two plane lattice walks both of which begin at the origin, which end at the same point, and which take $n$ steps each of which is N or E, and which intersect each other in exactly $k$ points. We will count such pairs.

From the proposition of section 1 above, we have that

$$(x + y + 2f)^{k+1} = \sum_{i,j} N_k^{i+j,j} x^i y^j.$$

If we put $y = x$ we get

$$2^{k+1}(x + f(x,x))^{k+1} = \sum_{i,j} N_k^{i+j,j} x^{i+j}$$

$$= (1 - \sqrt{1-4x})^{k+1} = (2x)^{k+1} \sum_{m \geq 0} \frac{(k+1)(2m+k)!}{m!(m+k+1)!} x^m.$$

Now matching $[x^n]$,

$$\sum_{i+j=n} N_k^{i+j,j} = 2^{k+1}(k+1) \frac{(2n-k-2)!}{n!(n-k-1)!}.$$

Since there are $\binom{2n}{n}$ pairs of walks that start at the origin and end at the same point, we can say this in probability language, like this:

*If two independent walkers start at the origin, and each takes $n$ steps, each step being N or E with equal probability, and if they finish at the same point, then the probability that the interior of their paths intersect in $k$ points is*

$$p(n,k) = \frac{2^{k+1}(k+1)(2n-k-2)!n!}{(n-k-1)!(2n)!} \qquad (n \geq 1; 0 \leq k \leq n-1).$$

The fact that $P(n) \stackrel{\text{def}}{=} \sum_{0 \leq k \leq n} p(n,k)$ is 1 for all $n \geq 1$ can be proved by adapting the WZ method (see [5]) as follows. (Some adaptation is necessary since the interval of summation does not coincide with the full support of the summand.) Begin with the routinely verifiable fact that

$$p(n+1, k) - p(n, k) = g(n, k+1) - g(n, k),$$

where $g(n,k) = -(k+2)p(n,k)/(2n+1)$, and sum from $k = 0$ to $n+1$. The right side telescopes to zero, and we have $P(n+1) - P(n) = 0$ for $n \geq 1$. ∎

It is interesting to note that $p(n, 1) = 2p(n, 0)$ if $n > 1$.

The following result does away with the condition that the walks finish at the same point.

**Theorem 3.** *The number of pairs of walks that begin at the origin and proceed with N or E steps, which each have n steps, and which intersect each other exactly k times excluding the origin, is $2^k \binom{2n-k}{n}$.*

**Corollary 1.** *Since there are $4^n$ pairs of walks of n steps, we have a combinatorial interpretation of the well known identity*

$$\sum_{k \geq 0} 2^k \binom{2n-k}{n} = 4^n.$$

**Corollary 2.** *The probability that two such walks do not intersect at all is $\binom{2n}{n}/4^n$.*

**Corollary 3.** *The average number of times that two independent random walks of n steps, beginning at the origin, cross each other is*

$$\frac{(2n+1)!}{4^n n!^2} - 1 = 2\sqrt{\frac{n}{\pi}} - 1 + o(1).$$

**Proof.** If $f(n,k)$ is the number of pairs of walks from $(0,0)$, of $n$ steps, that intersect $k$ times, and if $\phi(n,k)$ is the number of pairs of such walks that intersect $k$ times and finish at the same point, then

$$\sum_j \phi(j, k-1) f(n-j, 0) = f(n,k). \tag{8}$$

If we sum over $k$ we find that $4^n = \sum_j f(n-j, 0) \binom{2j}{j}$. So the generating function for pairs of walks with 0 intersections is $F_0(x) = (1-4x)^{-\frac{1}{2}}$, and $f(n,0) = \binom{2n}{n}$. Now from (8), since $\Phi_k(x) = (1 - \sqrt{1-4x})^{k+1}$, we have

$$F_k(x) = \frac{(1-\sqrt{1-4x})^k}{\sqrt{1-4x}} = 2^k \sum_j \binom{2j-k}{j} x^j,$$

hence $f(n,k) = 2^k \binom{2n-k}{n}$. ∎

For another proof one can look at the difference of the two walks. If $(x', y')$ and $(x'', y'')$ are the coordinates of the walkers on the two walks, put $(x, y) := (x' - x'', y' - y'')$. Then $(x,y)$ walks to $(x,y)$ with probability 1/2, to $(x+1, y-1)$ with probability 1/4, and to $(x-1, y+1)$ with probability 1/4. Thus the difference walk takes place entirely on the line $x + y = 0$. The statistics of intersections of the original pair of walks are identical with those of returns to 0 of a single one-dimensional walk of twice as many steps, and are well known.

### 6. Nonintersecting walks with a barrier

In this section we consider a variation on the original problem. For convenience, we will think of the walkers as moving to the origin instead of away from it.

So consider two walkers, $U$ and $L$, that start at the respective lattice points $(a, b+x+1)$ and $(a+x+1, b)$ (where $a, b, x \geq 0$). They independently move either West or South until they reach the $x$- or $y$-axis, where they are constrained to move along the axis to the origin. At each lattice point of $A = \{(r,s) : r, s > 0\}$ the probability of moving West is $p(r,s)$ and the probability of moving South is $1 - p(r,s)$. We want to know the probability $B(a,b,x)$ that the first time the walkers meet is at the origin. We say such a pair of walks is "valid."

Let us stop the walks after $a + b + x$ steps. Then a walker ends at either the point $(0, 1)$ or the point $(1, 0)$. The condition that a pair of walks is valid is thus equivalent to all three of the following holding: (1) the walkers never meet in $A$, (2) $U$ ends at $(0, 1)$, and (3) $L$ ends at $(1, 0)$.

Now, we claim that the probability that (2) and (3) are true but the walks intersect is equal to the probability that $U$ ends at $(1, 0)$ and $L$ ends at $(0, 1)$. For, if we have a pair $P$ of walks which satisfy (2) and (3) but intersect, we can create a new pair of walks $P'$ by interchanging the segments from the beginning to the first intersection ("initial segments") of the two walks in $P$. The resultant $P'$ is a pair of walks in which $U$ ends at $(1, 0)$ and $L$ ends at $(0, 1)$, where $P'$ has the same probability of occurrence as does $P$. Similarly, a pair $P'$ of walks in which $U$ ends at $(1, 0)$ and $L$ ends at $(0, 1)$ must intersect somewhere; if we interchange their initial segments we obtain a pair $P$ of walks in which $U$ ends at $(0, 1)$ and $L$ ends at $(1, 0)$.

Hence if $u$ denotes the probability that $U$ ends at $(0, 1)$, and $l$ denotes the probability that L ends at $(1, 0)$, then the probability that a pair of walks is valid is

$$B = ul - (1-u)(1-l) = u + l - 1.$$

Now, to simplify the above we proceed as follows. We remove the constraints on the axes and extend $p(r,s)$ to the points where $r \leq 0$ or $s \leq 0$ arbitrarily. Then the probability $u$ is the probability that a single walker starting at $(a, b+x+1)$ and moving either to the West, with probability $p(r,s)$, or South, with probability $1 - p(r,s)$, is, after $a+b+x$ steps, at one of the vertices in the set of lattice points on the line $x + y = 1$ that lie on or to the left of the $y$-axis. Similarly, $l$ is the probability that a single walker starting at $(a+x+1, b)$ is, after $a+b+x$ steps, in the set of lattice points on the line $x+y = 1$ that lie on or below the $x$-axis. If the transitional probabilities $p(r,s)$ in $A$ are a function of $r + s$ only, say $p(r,s) = p_{r+s}$, then we can extend $p(r,s)$ so this remains true. In this case, $l$ is the probability that a single walker starting at $(a, b+x+1)$ is, after $a+b+x$ steps, in the set of lattice points on the line $x+y = 1$ that lie on or below the point $(-x, 1+x)$. The expression $u + l - 1$ thus simplifies to yield the following:

**Theorem 4.** *If the transitional probability at $(r,s)$ is a function of $r + s$ only, then $B(a,b,x)$ is equal to the probability that a single unconstrained walker starting at $(a, b+x+1)$ and walking for $a+b+x$ steps ends up at one of the $x+1$ points $\{(-t, 1+t) : 0 \leq t \leq x\}$.*

**Corollary.** *If all the transitional probabilities have the same value $p$, then*

$$B(a,b,x) = \sum_{t=0}^{x} \binom{a+b+x}{a+t} p^{a+t} q^{b+x-t}. \qquad (q = 1-p)$$

In particular, the special case that originally interested us was the case of two walkers starting at the same point $(a+1, b+1)$ and moving as above: in this the answer is simply $2\binom{a+b}{a}p^{a+1}q^{b+1}$.

Brandeis University, Waltham, MA 02254-9110
University of Pennsylvania, Philadelphia, PA 19104-6395
New York Life Insurance Co., New York, NY 10010
University of Pennsylvania, Philadelphia, PA 19104-6395
University of Waterloo, Waterloo, Ontario, Canada N2L 3G1